\documentclass[twoside,10pt]{article}

\newcommand{\qed}{\quad \rule[-0.5mm]{1.9mm}{3.2mm}}
\newcommand{\equsep}{\rule{0mm}{6mm}}

\usepackage{fancyheadings}

\newtheorem{thm}{Theorem}[section]
\newtheorem{prop}[thm]{Proposition}
\newtheorem{cor}[thm]{Corrolary}

\newtheorem{lem}[thm]{Lemma}
\newtheorem{conj}[thm]{Conjecture}
\newtheorem{problem}[thm]{Problem}

\def\authorname{W. Y. C. Chen and V. J. W. Guo}
\def\papername{Bijections behind the Ramanujan Polynomials}
\newcommand{\newsection}[1]{\section{#1} \setcounter{equation}{0}}

\begin{document}
\thispagestyle{plain}

\begin{center}
{\Large\bf Bijections behind the Ramanujan Polynomials}
\vskip 3mm
Dedicated to Professor Dominique Foata, in honor of his 65th  Birthday
\end{center}

\noindent
William Y. C. Chen\\
Center for Combinatorics\\
The Key Laboratory of Pure Mathematics and
    Combinatorics  of Ministry of Education\\
Nankai University, Tianjin 300071, P. R. China\\
Email: chenstation@hotmail.com

\vskip 6pt
\noindent
and

\vskip 3mm
\noindent
Theoretical Division, MS K710\\
Los Alamos National Laboratory\\
Los Alamos, New Mexico 87545, USA

\vskip 6pt
\noindent
and

\vskip 3mm \noindent {Victor J. W. Guo}
\\ Center for
Combinatorics\\ The Key Laboratory of Pure Mathematics and
    Combinatorics  of Ministry of Education\\
Nankai University, Tianjin 300071, P. R. China\\
Email: jwguo@eyou.com

\vskip 5mm

\begin{minipage}{5in}
 {\bf Abstract.}
The Ramanujan polynomials were introduced by Ramanujan in his
study of power series inversions. 
In an approach to the Cayley formula on the number of trees, Shor discovers
a refined recurrence relation in terms of the number of improper edges,
without realizing the connection to the Ramanujan polynomials. On
the other hand, Dumont and Ramamonjisoa independently take the
grammatical approach to a sequence associated with the Ramanujan
polynomials and have reached the same conclusion as Shor's.
It was a 
coincidence for Zeng to realize that the Shor polynomials turn out
to be the Ramanujan polynomials through an explicit substitution
of parameters. 
Shor also discovers a recursion of Ramanujan
polynomials which is equivalent to the Berndt-Evans-Wilson
recursion under the substitution of Zeng, and asks for a
combinatorial interpretation. The objective of this paper is to
present a bijection for the Shor recursion, or and
Berndt-Evans-Wilson recursion, answering the question of Shor.
Such a bijection also leads to  a combinatorial interpretation of
the recurrence relation originally given by Ramanujan.

\end{minipage}

\newsection{Introduction}

The original Ramanujan polynomials $\psi_k(r,x)$, where $r$ is
any nonnegative integer and $1\leq k \leq r+1$,  are
defined by the following generating function equation:
\begin{equation}
\sum\limits_{k=0}^{\infty}\frac{(x+k)^{r+k}e^{-u(x+k)}u^k}{k!}
=\sum\limits_{k=1}^{r+1}\frac{{\psi}_k(r,x)}{(1-u)^{r+k}}.
\label{def-gen-id}
\end{equation}
Ramanujan gives a recurrence relation
of $\psi_{k}(r,x)$ as follows:
\begin{equation}
{\psi}_k(r+1,x)=
(x-1){\psi}_k(r,x-1)+{\psi}_{k-1}(r+1,x)-{\psi}_{k-1}(r+1,x-1).
\label{1-rec-org}
\end{equation}
where $1 \leq k \leq r+1, \,  {\psi}_{1}(0,x)=1,$ and
$ {\psi}_{k}(r,x)=0$  if $k \not \in [r+1]$.
Note that here we have adopted the standard notation
$[n]:=\{1, 2, \ldots, n\}$ for a positive integer $n$.

Berndt et al. \cite{Berndt83,Berndt85} find an
elegant proof of (\ref{def-gen-id}) justifying
the existence of the polynomials
$\psi_k(r,x)$ and obtain the following recurrence relation:
\begin{equation}
 {\psi}_{k}(r,x)=(x-r-k+1){\psi}_{k}(r-1,x)+(r+k-2){\psi}_{k-1}(r-1,x),
\label{2-rec-org}
\end{equation}
where the initial value of $\psi_k(r,x)$ and the
ranges of indices are given as above.

It is worth noting that the Ramanujan polynomials satisfy the
following identity:
\begin{equation}
\sum_{k=1}^{r+1} \; \psi_{k}(r, x) = x^r.
\label{ram-xr}
\end{equation}

\[ \mbox{ Table of $\psi_{k}(r, x)$. }\]
\centerline{
\footnotesize{
\begin{tabular}{|l|r|r|r|r|r|}\hline
     $ k\backslash r$ & $0$ & $1$   & $2$ & $3$ & $4$\\ \hline
      $1$             & $1$ & $x-1$ & $x^2-3x+2$ & $x^3-6x^2+11x-6$
                    & $x^4-10x^3+35x^2-50x+24$       \\ \hline
      $2$             &     & $1$   & $3x-5$ & $6x^2-26x+26$
                    & $10x^3-80x^2+200x-154$          \\ \hline
      $3$             &     &       & $3$    & $15x-35$
                    & $45x^2-255x+340$               \\ \hline
      $4$             &     &       &        & $15$
                    & $105x-315$                     \\ \hline
      $5$             &     &       &        &
                    & $105$                          \\ \hline
     $ \sum_k$        & $1$ & $x$   & $x^2$  & $x^3$
                    & $x^4$                          \\ \hline
   \end{tabular}
             }
           }
\vskip 0.5cm

It turns out that the Ramanujan polynomials coincide with the
polynomials $Q_{n,k}(x)$ introduced by Shor [\ref{Shor95}], where
$n\geq 1$, and $0\leq  k \leq n-1$. Moreover, for $n=0$
or $k\not \in [n-1]$, we define $Q_{n,k}(x)$ to be zero. Shor's
recursive definition of $Q_{n,k}(x)$ goes as follows:
\begin{equation}
Q_{n,k}(x)=(x+n-1)Q_{n-1,k}(x)+(n+k-2)Q_{n-1,k-1}(x),
\label{2-rec-new}
\end{equation}
for $n\geq 1$ and $k\leq n-1$, where $Q_{1,0}(x) =1$ and
$Q_{n,k}(x)=0$ if $k\geq n$ or $k<0$.
Zeng [\ref{Zeng}, Proposition 7] establishes the following remarkable connection:
\begin{equation}
Q_{n,k}(x) = \psi_{k+1}(n-1, x+n).
\label{zeng-sub}
\end{equation}
The tree enumeration flavor of $Q_{n,k}(x)$ is evidenced by the following identity:
\begin{equation}
\sum_{k=0}^{n-1}\, Q_{n,k}(x) = (x+n)^{n-1}.
\label{qnkx-sum}
\end{equation}

\[ \mbox{ Table of $Q_{n,k}(x)$. }\]
\centerline{
\footnotesize{
   \begin{tabular}{|l|r|r|r|r|r|}\hline
     $ k\backslash n$ & $1$ & $2$   & $3$ & $4$ & $5$\\ \hline
      $0$             & $1$ & $x+1$ & $x^2+3x+2$ & $x^3+6x^2+11x+6$
                    & $x^4+10x^3+35x^2+50x+24$       \\ \hline
      $1$             &     & $1$   & $3x+4$ & $6x^2+22x+18$
                    & $10x^3+70x^2+150x+96$          \\ \hline
      $2$             &     &       & $3$    & $15x+25$
                    & $45x^2+195x+190$               \\ \hline
      $3$             &     &       &        & $15$
                    & $105x+210$                     \\ \hline
      $4$             &     &       &        &
                    & $105$                          \\ \hline
     $ \sum_k$ & $1$ & $x+2$   & $(x+3)^2$ & $(x+4)^3$
                    & $(x+5)^4$                      \\ \hline
   \end{tabular}
             }
           }
\vskip 0.5cm

In his approach to the enumeration of trees, Shor [\ref{Shor95}] has considered  the following recurrence relation:
\begin{equation}
      f_{n,k}=(n-1)f_{n-1,k}+(n+k-2)f_{n-1,k-1},
\label{shor-rec-n}
\end{equation}
where $f_{1,0}=1$, $n\geq 1$, $k\leq n-1$, and $f_{n,k}=0$
otherwise. One sees that $f_{n,k}$ is
the value of $Q_{n,k}(x)$ evaluated at $x=0$, and
that $f_{n,k}$ satisfies the following identity:
\begin{equation}
\sum_{k=0}^{n-1} \; f_{n,k} = n^{n-1}.
\label{fnk-sum}
\end{equation}
Shor shows that
$f_{n,k}$  is in fact the number of rooted trees
on $[n]$ with $k$ improper edges.
However, he did not seem to have
noticed the connection of his formula to
the work of Ramanujan.
On the other hand,
Dumont and Ramamonjisoa \cite{Dumont} use the grammatical
method introduced by Chen in \cite{Chen93} to obtain the same
combinatorial interpretation.

Besides the recurrence relation (\ref{2-rec-new}) for
$Q_{n,k}(x)$, Shor \cite{Shor95} derives the following recurrence
relation, and asks for a combinatorial interpretation:
\begin{equation}
Q_{n,k}(x)=(x-k+1)Q_{n-1,k}(x+1)+(n+k-2)Q_{n-1,k-1}(x+1).
\label{shor-rec}
\end{equation}
The above recurrence relation turns out to be
equivalent to the Berndt-Evans-Wilson recursion (\ref{2-rec-org})
by the substitution (\ref{zeng-sub}) of Zeng.

The aim of this paper is to
 construct a bijection for (\ref{shor-rec}), answering the question of Shor.
We note that the above relation is indeed the
same as the recurrence relation (\ref{1-rec-org})
under the substitution (\ref{zeng-sub}) of Zeng.
Therefore,  we also obtain a combinatorial interpretation
of the recurrence relation (\ref{1-rec-org})
originally presented by Ramanujan.

\newsection{The Zeng Interpretations and the Shor Recursion}

We will follow most notation in Zeng \cite{Zeng}. The set of
rooted labeled trees on $[n]$ is denoted by $\mathcal{R}_n$ . If
$T\in \mathcal{R}_n$, and $x$ is a node of $T$, the subtree rooted
at $x$ is denoted by $T_x$. We let $\beta(x)$, or $\beta_{T}(x)$
be the smallest node on $T_x$. For notational simplicity, we also
use $\beta_T$ or $\beta(T)$ to denote the minimum element in $T$,
and we  sometimes write $T(x)$ for $T_x$ in the purpose of
avoiding multiple subscripts. We  say that a node $z$ of $T$ is a
descendant of $x$, (or $x$ is an ancestor of $z$), if $z$ is a
node of $T_x$. In particular, each node is a descendant of itself.
For any edge $e=(x,y)$ of a tree $T$, if $y$ is a node of $T_x$,
we call $x$ the {\it father node} of $e$, $y$ the {\it child node}
of $e$, $x$ the father of $y$, and $y$ a child of $x$. Assume
$e=(x,y)$ is an edge of a tree $T$, and $y$ is the child node of
$e$, we say that $e$ is a {\it proper} edge, if $x<\beta_{T}(y)$.
Otherwise, we call $e$ an {\it improper} edge. The degree of a
node $x$ in a rooted tree $T$ is the number of children of $x$,
and is denoted by $\deg(x)$, or $\deg_{T}(x)$.  An unrooted
labeled tree will be treated as a rooted tree in which the
smallest node is chosen as the root. Then the above definitions
are still valid for unrooted trees. Denote by $\mathcal{T}_{n,k}$
and $\mathcal{R}_{n,k}$ the sets of labeled trees and rooted
labeled trees on $[n]$ with $k$ improper edges, respectively.
Moreover, we may impose some conditions on the sets
$\mathcal{T}_{n,k}$ and $\mathcal{R}_{n,k}$ to denote the subsets
of trees that satisfy these conditions. For example,
$\mathcal{T}_{n+1,k} [\deg(n+1)=0]$ stands for the subset of
$\mathcal{T}_{n+1,k}$ subject to the condition $\deg(n+1)=0$.

\setlength{\unitlength}{0.08in}
\begin{figure}[ht]
\begin{picture}(60,10)(10,0)

    \put(41,6){\circle*{0.7}}       \put(40.1,6.8){\makebox(2,1)[l]{14}}
    \put(38,2){\circle*{0.7}}       \put(39.3,1.9){\makebox(2,1)[l]{6}}
    \put(49,2){\circle*{0.7}}       \put(48.6,0){\makebox(2,1)[l]{13}}
    \put(30,-2){\circle*{0.7}}       \put(29.7,-4){\makebox(2,1)[l]{5}}
    \put(35,-2){\circle*{0.7}}       \put(34.7,-4){\makebox(2,1)[l]{3}}
    \put(38,-2){\circle*{0.7}}       \put(38.6,-2.2){\makebox(2,1)[l]{12}}
    \put(42,-2){\circle*{0.7}}       \put(41.6,-4){\makebox(2,1)[l]{4}}
    \put(46,-2){\circle*{0.7}}       \put(45.6,-1.2){\makebox(2,1)[l]{8}}
    \put(32,-6){\circle*{0.7}}       \put(31.7,-8){\makebox(2,1)[l]{10}}
    \put(36,-6){\circle*{0.7}}       \put(36,-8){\makebox(2,1)[l]{16}}
    \put(40,-6){\circle*{0.7}}       \put(40.7,-6.2){\makebox(2,1)[l]{19}}
    \put(46,-6){\circle*{0.7}}       \put(45.5,-8){\makebox(2,1)[l]{7}}
    \put(50,-6){\circle*{0.7}}       \put(49.7,-8){\makebox(2,1)[l]{2}}
    \put(54,-6){\circle*{0.7}}       \put(53.7,-5.2){\makebox(2,1)[l]{17}}
    \put(29,-10){\circle*{0.7}}      \put(28.7,-12){\makebox(2,1)[l]{20}}
    \put(34,-10){\circle*{0.7}}       \put(33.1,-12){\makebox(2,1)[l]{18}}
    \put(42,-10){\circle*{0.7}}       \put(41.7,-12){\makebox(2,1)[l]{1}}
    \put(62,-10){\circle*{0.7}}       \put(61.7,-12){\makebox(2,1)[l]{9}}
    \put(50,-10){\circle*{0.7}}       \put(49.7,-12){\makebox(2,1)[l]{15}}
    \put(26,-14){\circle*{0.7}}       \put(25.1,-16){\makebox(2,1)[l]{11}}

  \put(40.8,6){\line(-3,-4){6}}   \put(41.2,6){\line(-3,-4){6}}
  \put(28.8,-10){\line(-3,-4){3}} \put(29.2,-10){\line(-3,-4){3}}
  \put(40.7,6){\line(2,-1){8}}    \put(41.3,6){\line(2,-1){8}}
  \put(37.7,2){\line(-2,-1){8}}   \put(38.3,2){\line(-2,-1){8}}
  \put(37.85,2){\line(0,-1){4}}   \put(38.15,2){\line(0,-1){4}}
  \put(35,-2){\line(-3,-4){6}}
  \put(37.8,2){\line(1,-1){4}}    \put(38.2,2){\line(1,-1){4}}
  \put(42,-2){\line(1,-1){8}}
  \multiput(37.7,2)(16,-8){2}{\line(2,-1){8}}
  \multiput(38.3,2)(16,-8){2}{\line(2,-1){8}}
  \put(46,-2){\line(2,-1){8}}
  \put(38,-2){\line(-1,-2){4}}
  \put(37.8,-2){\line(1,-2){4}}   \put(38.2,-2){\line(1,-2){4}}
  \put(45.8,-2){\line(1,-1){4}}   \put(46.2,-2){\line(1,-1){4}}
 \end{picture}
 \vskip 3.2cm
 \caption{Improper Edges Shown as Double Edges}
\end{figure}

\begin{thm}{\rm \bf(Zeng [6,\,Propositions 1,\,2,\,7])}
The polynomials $Q_{n,k}(x)$ have the following interpretations:
\begin{eqnarray}
Q_{n,k}(x)&=&
\sum\limits_{T\in\mathcal{T}_{n+1,k}}x^{\deg_{T}(1)-1}.
\label{zeng-1}   \\
Q_{n,k}(x)
&=&\sum\limits_{T\in\mathcal{R}_{n,k}}(x+1)^{\deg_{T}(1)}.
\label{zeng-2}
\end{eqnarray}
\end{thm}

In fact, the above theorem can be reformulated by the following
relations:
\begin{eqnarray}
(x+n-1)Q_{n-1,k}(x) & = &
\sum\limits_{T\in\mathcal{T}_{n+1,k}
[\deg(n+1)=0]}x^{\deg_{T}(1)-1},
\label{q-1} \\
(n+k-2)Q_{n-1,k-1}(x) & = &
\sum\limits_{T\in\mathcal{T}_{n+1,k}
[\deg(n+1)>0]} \, x^{\deg_{T}(1)-1} .
\label{q-2}
\end{eqnarray}

Zeng \cite{Zeng} proves the two interpretations  (\ref{zeng-1})
and (\ref{zeng-2}) of $Q_{n,k}(x)$ by similar arguments. One
naturally expects to make a combinatorial connection bridging
these two formulations, and this consideration was mentioned by
Zeng. We now provide such an argument for the equivalence between
(\ref{zeng-1}) and (\ref{zeng-2}), that is,
\begin{equation}
\sum\limits_{T\in\mathcal{T}_{n+1,k}}x^{\deg_T(1)-1}
 =\sum\limits_{T\in\mathcal{R}_{n,k}}(x+1)^{\deg_T(1)}.
\label{equiv-12}
\end{equation}

{\em Proof.}
Let us consider the binomial expansion of the right hand side
of (\ref{equiv-12}).  The binomial expansion can
be visualized by coloring
the children of the node $1$ with
black and white colors. Let $T$ be a rooted tree in $\mathcal{R}_{n,k}$,
and let $T$ have the children of $1$ colored in either black or white. Let $B$ be the set of
children of $1$ in $T$ which are colored in black.
Now we may introduce a new node $0$, and move
the subtrees of $T$ rooted at the nodes in $B$ as
the subtrees of $0$, and moreover, move the remaining
subtree of $T$ as a subtree of $0$.  Therefore, we obtain a
rooted tree on $\{ 0, 1, 2, \ldots, n\}$, say $T'$.
Note that the children of $0$ which come from the
black nodes can be
easily distinguished from the child of $0$ which
is the original root of $T$ because the node $1$
remains in the subtree of original root.
Finally, if we relabel the set $\{0, 1, 2, \ldots, n\}$
by the set $[n+1]$, namely, relabeling $i$ by $i+1$,
we get an unrooted tree on $[n+1]$ which preserves the
number of improper edges.
Furthermore, one sees that the above construction can
be reversed. This completes the proof. \qed

\begin{cor} We have
\begin{equation}
 Q_{n,k}(x-1)=\sum\limits_{T\in\mathcal{T}_{n+1,k}
 [\deg(2)=0] } \, x^{\deg_T(1)-1}.
\end{equation}
\end{cor}

{\em Proof.} It follows from $(\ref{zeng-2})$ that
\begin{equation}
Q_{n,k}(x-1)=\sum\limits_{T\in\mathcal{R}_{n,k}}\,
x^{\deg_T(1)}.
\end{equation}
We now construct a bijection from $\mathcal{R}_{n,k}[\deg(1)=r] $
to $\mathcal{T}_{n+1,k}[\deg(1)=r+1, \deg(2)=0]$. Given
$T\in\mathcal{R}_{n,k}[\deg(1)=r]$, we now introduce a new root
$0$, and put $T$ as a subtree of $0$. Then we move all the
subtrees of $1$ and make them as subtrees of $0$. Finally, by
relabeling a node $i$ by $i+1$, we obtain a tree
$T'\in\mathcal{T}_{n+1,k}[\deg(1)=r+1,
\deg(2)=0]$. It is clear that the construction is reversible. This
completes the proof. \qed

Substituting $k$ by $k+1$ in (\ref{q-2}), we obtain
\begin{equation}
(n+k-1)Q_{n-1,k}(x)
= \sum\limits_{T\in\mathcal{T}_{n+1,k+1}
[\deg(n+1)>0]}\, x^{\deg_{T}(1)-1} .
\label{q-3}
\end{equation}

We are now ready to give another  combinatorial formulation of
the Shor recurrence relation (\ref{shor-rec}). Rewriting
(\ref{shor-rec}), by substituting $x$ with $x-1$, we get:
\begin{equation}
Q_{n,k}(x-1) = (x-k) Q_{n-1, k}(x) + (n+k-2) Q_{n-1, k-1}(x).
\label{q-4}
\end{equation}

If we express the term $(x-k) Q_{n-1, k}(x)$
as
\[ (x+n-1) Q_{n-1, k}(x) - [n +(k+1) -2] Q_{n-1, (k+1)-1}(x)\, ,\]
then the Shor recurrence relation (\ref{shor-rec})
is equivalent to the following combinatorial identity.

\begin{thm}\label{thm-main-comb}
For $n\geq 1$, and $0\leq k \leq n-1$,
we have
\begin{equation}
  \sum\limits_{T\in\mathcal{T}_{n+1,k}
[\deg(2)>0]} \, x^{\deg_T(1)-1}
 =\sum\limits_{T\in\mathcal{T}_{n+1,k+1}
[\deg(n+1)>0]}\, x^{\deg_T(1)-1}.
\label{main-comb}
\end{equation}
\end{thm}

We now present an inductive proof of the above fact, while the
next section will be engaged in a purely combinatorial treatment.
Clearly, for $n\geq 1$,
(\ref{main-comb}) can be restated as follows with the
notation $T_{n,k}[\cdots]:=|\mathcal{T}_{n,k}[\cdots]|$:
\begin{equation}
 T_{n+1,k}[\deg(2)>0, \deg(1)=r] =
 T_{n+1,k+1}[\deg(n+1)>0, \deg(1)=r] . \label{tnk-r}
\end{equation}

{\em Proof.}
For $n\geq 2$, the arguments of
Shor \cite{Shor95} or Zeng \cite{Zeng} imply the following
identities:
\begin{eqnarray*}
\mbox{(i)}&&
T_{n+1,k+1}[\deg(n+1)>0,\deg(1)=r]=(n+k-1)T_{n,k}[\deg(1)=r].\\[5pt]
\mbox{(ii)}&&T_{n+1,k}[\deg(2)>0,\deg(1)=r]\\[5pt]
&&\quad =\,
  (n-2)T_{n,k}[\deg(2)>0,\deg(1)=r]+T_{n,k}[\deg(2)>0,\deg(1)=r-1]\\[5pt]
&&\quad \qquad  +\,
T_{n,k}[\deg(1)=r]+(n+k-2)T_{n,k-1}[\deg(2)>0,\deg(1)=r]. \\[5pt]
\mbox{(iii)}&&T_{n+1,k}[\deg(1)=r]\\[5pt] &&\quad =\,
(n-1)T_{n,k}[\deg(1)=r]+T_{n,k}[\deg(1)=r-1] \\[5pt] &&\quad
\qquad+ \, (n+k-2)T_{n,k-1}[\deg(1)=r].
\end{eqnarray*}
Because of (i), (\ref{tnk-r}) can be deduced from the following
relation:
\begin{equation}
T_{n+1,k}[\deg(2)>0,\deg(1)=r]=(n+k-1)T_{n,k}[\deg(1)=r],
\label{deg-2-0}
\end{equation}
for $n\geq 1$. The above claimed identity obviously holds for $n=1$.
 Suppose (\ref{deg-2-0}) holds
for $n-1$. From (i) -- (iii) and the inductive hypothesis,
it follows that
\begin{eqnarray*}
\lefteqn{T_{n+1,k}[\deg(2)>0,\deg(1)=r]}\\[5pt]
&=&(n-2)T_{n,k}[\deg(2)>0,\deg(1)=r]+T_{n,k}[\deg(2)>0,\deg(1)=r-1]\\[5pt]
&&\quad + \,
T_{n,k}[\deg(1)=r]+(n+k-2)T_{n,k-1}[\deg(2)>0,\deg(1)=r]\\[5pt]
&=&
(n-2)(n+k-2)T_{n-1,k}[\deg(1)=r]+(n+k-2)T_{n-1,k}[\deg(1)=r-1]\\[5pt]
&&\quad +\,
T_{n,k}[\deg(1)=r]+(n+k-2)(n+k-3)T_{n-1,k-1}[\deg(1)=r]\\[5pt] &=&
(n+k-2)
\bigl\{(n-2)T_{n-1,k}[\deg(1)=r]+T_{n-1,k}[\deg(1)=r-1]\\[5pt]
&&\quad +(n+k-3)T_{n-1,k-1}[\deg(1)=r] \bigr\}
+T_{n,k}[\deg(1)=r]\\[5pt]
&=&(n+k-2)T_{n,k}[\deg(1)=r]+T_{n,k}[\deg(1)=r]\\[5pt]
&=&(n+k-1)T_{n,k}[\deg(1)=r].
\end{eqnarray*}
Thus (\ref{deg-2-0}) holds for $n$.
This completes the proof. \qed

We further remark that the following recurrence relations
presented by Zeng \cite{Zeng} also follow from
the above combinatorial identity:
\begin{eqnarray}
Q_{n,k}(x)& =& (x+n-1)Q_{n-1,k}(x)+Q_{n,k-1}(x)-Q_{n,k-1}(x-1),
\label{q-f-1}
\\
\equsep
Q_{n,k}(x)& =& Q_{n,k}(x-1)+(n+k-1)Q_{n-1,k}(x).
\end{eqnarray}

Note that the recurrence relation (\ref{q-f-1}) is equivalent to
the original Ramanujan recursion (\ref{1-rec-org}). A bijective
proof of (\ref{main-comb}) will be the objective of the next
section.

\newsection{The Bijections}

In order to demonstrate (\ref{main-comb}) combinatorially, it
would be ideal to directly construct a bijection from
$\mathcal{T}_{n+1,k} [\deg(2)>0]$ to $\mathcal{T}_{n+1,k+1}
[\deg(n+1)>0]$ which preserves the degree of $1$. Although it
looks that such a bijection should be easy to construct by moving
a child of $2$ in a tree in $\mathcal{T}_{n+1,k} [\deg(2)>0]$ to
the node $n+1$, achieving such a task turns out to be quite subtle. To
achieve this goal, we first find a stronger bijection on rooted
trees subject to certain degree constraints.

\begin{thm}\label{thm-bij-r}
For $n\geq 1$ and $0 \leq k < n$, we have
the following bijection:
\begin{equation}
      \mathcal{R}_{n,k}[\deg(1)>0]
     \longleftrightarrow
      \mathcal{R}_{n,k+1}[\deg(n)>0].
\label{bij-r}
\end{equation}
\end{thm}

Here is an example for $n=4$, $k=1$. There are $16$ trees for each
side of (\ref{bij-r}). The trees in $\mathcal{R}_{4,1}[\deg(1)>0]$
are listed in Figure \ref{fig-deg1}.

\setlength{\unitlength}{0.08in}
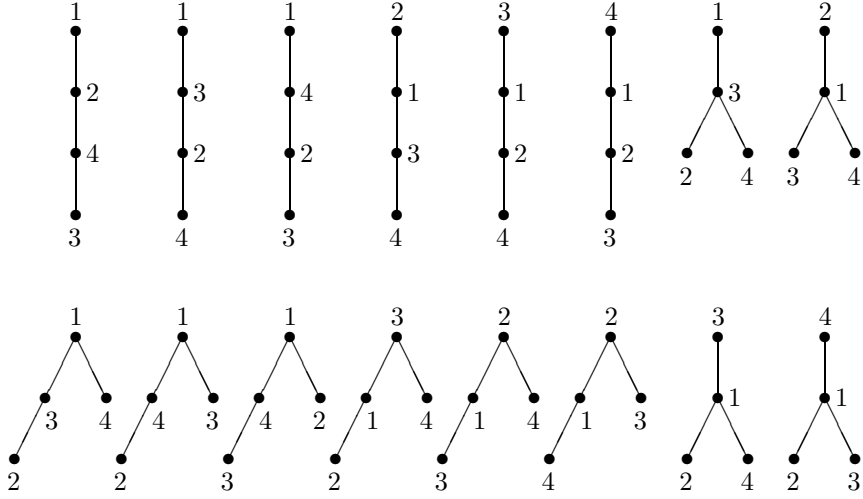
\begin{figure}[ht]
\hskip 1cm
\begin{picture}(60,20)

 \multiput(4,16)(7,0){6}{\circle*{0.7}} \put(3.6,16.8){\makebox(2,1)[l]{1}}
 \multiput(4,12)(7,0){6}{\circle*{0.7}} \put(4.7,11.5){\makebox(2,1)[l]{2}}
 \multiput(4,8)(7,0){6}{\circle*{0.7}}  \put(4.7,7.5){\makebox(2,1)[l]{4}}
 \multiput(4,4)(7,0){6}{\circle*{0.7}}  \put(3.5,2){\makebox(2,1)[l]{3}}
 \multiput(4,16)(7,0){6}{\line(0,-1){12}}

 \put(10.6,16.8){\makebox(2,1)[l]{1}}
 \put(11.7,11.5){\makebox(2,1)[l]{3}}
 \put(11.7,7.5){\makebox(2,1)[l]{2}}
 \put(10.5,2){\makebox(2,1)[l]{4}}

 \put(17.6,16.8){\makebox(2,1)[l]{1}}
 \put(18.7,11.5){\makebox(2,1)[l]{4}}
 \put(18.7,7.5){\makebox(2,1)[l]{2}}
 \put(17.5,2){\makebox(2,1)[l]{3}}

 \put(24.6,16.8){\makebox(2,1)[l]{2}}
 \put(25.7,11.5){\makebox(2,1)[l]{1}}
 \put(25.7,7.5){\makebox(2,1)[l]{3}}
 \put(24.5,2){\makebox(2,1)[l]{4}}

 \put(31.6,16.8){\makebox(2,1)[l]{3}}
 \put(32.7,11.5){\makebox(2,1)[l]{1}}
 \put(32.7,7.5){\makebox(2,1)[l]{2}}
 \put(31.5,2){\makebox(2,1)[l]{4}}

 \put(38.6,16.8){\makebox(2,1)[l]{4}}
 \put(39.7,11.5){\makebox(2,1)[l]{1}}
 \put(39.7,7.5){\makebox(2,1)[l]{2}}
 \put(38.5,2){\makebox(2,1)[l]{3}}

 \multiput(4,-4)(7,0){6}{\circle*{0.7}}  \put(3.6,-3.2){\makebox(2,1)[l]{1}}
 \multiput(2,-8)(7,0){6}{\circle*{0.7}}  \put(2,-10){\makebox(2,1)[l]{3}}
 \multiput(0,-12)(7,0){6}{\circle*{0.7}} \put(-0.5,-14){\makebox(2,1)[l]{2}}
 \multiput(6,-8)(7,0){6}{\circle*{0.7}}  \put(5.5,-10){\makebox(2,1)[l]{4}}
 \multiput(4,-4)(7,0){6}{\line(-1,-2){4}}
 \multiput(4,-4)(7,0){6}{\line(1,-2){2}}

  \put(10.6,-3.2){\makebox(2,1)[l]{1}}
  \put(9,-10){\makebox(2,1)[l]{4}}
  \put(6.5,-14){\makebox(2,1)[l]{2}}
  \put(12.5,-10){\makebox(2,1)[l]{3}}

  \put(17.6,-3.2){\makebox(2,1)[l]{1}}
  \put(16,-10){\makebox(2,1)[l]{4}}
  \put(13.5,-14){\makebox(2,1)[l]{3}}
  \put(19.5,-10){\makebox(2,1)[l]{2}}

  \put(24.6,-3.2){\makebox(2,1)[l]{3}}
  \put(23,-10){\makebox(2,1)[l]{1}}
  \put(20.5,-14){\makebox(2,1)[l]{2}}
  \put(26.5,-10){\makebox(2,1)[l]{4}}

  \put(31.6,-3.2){\makebox(2,1)[l]{2}}
  \put(30,-10){\makebox(2,1)[l]{1}}
  \put(27.5,-14){\makebox(2,1)[l]{3}}
  \put(33.5,-10){\makebox(2,1)[l]{4}}

  \put(38.6,-3.2){\makebox(2,1)[l]{2}}
  \put(37,-10){\makebox(2,1)[l]{1}}
  \put(34.5,-14){\makebox(2,1)[l]{4}}
  \put(40.5,-10){\makebox(2,1)[l]{3}}

 \multiput(46,16)(7,0){2}{\circle*{0.7}}  \put(45.6,16.8){\makebox(2,1)[l]{1}}
 \multiput(46,12)(7,0){2}{\circle*{0.7}}  \put(46.7,11.5){\makebox(2,1)[l]{3}}
 \multiput(44,8)(7,0){2}{\circle*{0.7}}   \put(43.5,6){\makebox(2,1)[l]{2}}
 \multiput(48,8)(7,0){2}{\circle*{0.7}}   \put(47.5,6){\makebox(2,1)[l]{4}}
 \multiput(46,16)(7,0){2}{\line(0,-1){4}}
 \multiput(46,12)(7,0){2}{\line(-1,-2){2}}
 \multiput(46,12)(7,0){2}{\line(1,-2){2}}

 \multiput(46,-4)(7,0){2}{\circle*{0.7}}  \put(45.6,-3.2){\makebox(2,1)[l]{3}}
 \multiput(46,-8)(7,0){2}{\circle*{0.7}}  \put(46.7,-8.5){\makebox(2,1)[l]{1}}
 \multiput(44,-12)(7,0){2}{\circle*{0.7}} \put(43.5,-14){\makebox(2,1)[l]{2}}
 \multiput(48,-12)(7,0){2}{\circle*{0.7}} \put(47.5,-14){\makebox(2,1)[l]{4}}
 \multiput(46,-4)(7,0){2}{\line(0,-1){4}}
 \multiput(46,-8)(7,0){2}{\line(-1,-2){2}}
 \multiput(46,-8)(7,0){2}{\line(1,-2){2}}

  \put(52.6,16.8){\makebox(2,1)[l]{2}}
  \put(53.7,11.5){\makebox(2,1)[l]{1}}
  \put(50.5,6){\makebox(2,1)[l]{3}}
  \put(54.5,6){\makebox(2,1)[l]{4}}

  \put(52.6,-3.2){\makebox(2,1)[l]{4}}
  \put(53.7,-8.5){\makebox(2,1)[l]{1}}
  \put(50.5,-14){\makebox(2,1)[l]{2}}
  \put(54.5,-14){\makebox(2,1)[l]{3}}

 \end{picture}
 \vskip 3cm
 \caption{16 trees in $\mathcal{R}_{4,1}[\deg(1)>0]$}
 \label{fig-deg1}
 \end{figure}

The trees in  $\mathcal{R}_{4,2}[\deg(4)>0]$ are in Figure
\ref{fig-deg4}.

\setlength{\unitlength}{0.08in}
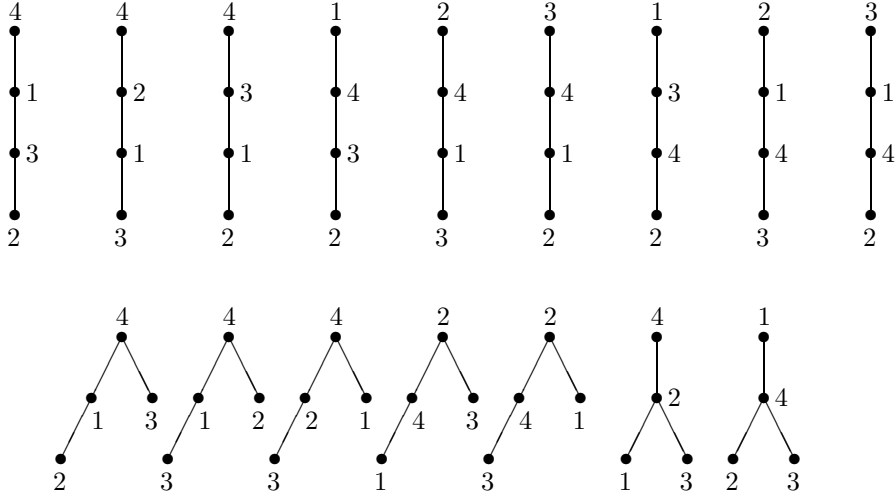
\begin{figure}[ht]
\hskip 0.3cm
\begin{picture}(60,20)

 \multiput(4,16)(7,0){9}{\circle*{0.7}} \put(3.6,16.8){\makebox(2,1)[l]{4}}
 \multiput(4,12)(7,0){9}{\circle*{0.7}} \put(4.7,11.5){\makebox(2,1)[l]{1}}
 \multiput(4,8)(7,0){9}{\circle*{0.7}}  \put(4.7,7.5){\makebox(2,1)[l]{3}}
 \multiput(4,4)(7,0){9}{\circle*{0.7}}  \put(3.5,2){\makebox(2,1)[l]{2}}
 \multiput(4,16)(7,0){9}{\line(0,-1){12}}

 \put(10.6,16.8){\makebox(2,1)[l]{4}}
 \put(11.7,11.5){\makebox(2,1)[l]{2}}
 \put(11.7,7.5){\makebox(2,1)[l]{1}}
 \put(10.5,2){\makebox(2,1)[l]{3}}

 \put(17.6,16.8){\makebox(2,1)[l]{4}}
 \put(18.7,11.5){\makebox(2,1)[l]{3}}
 \put(18.7,7.5){\makebox(2,1)[l]{1}}
 \put(17.5,2){\makebox(2,1)[l]{2}}

 \put(24.6,16.8){\makebox(2,1)[l]{1}}
 \put(25.7,11.5){\makebox(2,1)[l]{4}}
 \put(25.7,7.5){\makebox(2,1)[l]{3}}
 \put(24.5,2){\makebox(2,1)[l]{2}}

 \put(31.6,16.8){\makebox(2,1)[l]{2}}
 \put(32.7,11.5){\makebox(2,1)[l]{4}}
 \put(32.7,7.5){\makebox(2,1)[l]{1}}
 \put(31.5,2){\makebox(2,1)[l]{3}}

 \put(38.6,16.8){\makebox(2,1)[l]{3}}
 \put(39.7,11.5){\makebox(2,1)[l]{4}}
 \put(39.7,7.5){\makebox(2,1)[l]{1}}
 \put(38.5,2){\makebox(2,1)[l]{2}}

 \put(45.6,16.8){\makebox(2,1)[l]{1}}
 \put(46.7,11.5){\makebox(2,1)[l]{3}}
 \put(46.7,7.5){\makebox(2,1)[l]{4}}
 \put(45.5,2){\makebox(2,1)[l]{2}}

 \put(52.6,16.8){\makebox(2,1)[l]{2}}
 \put(53.7,11.5){\makebox(2,1)[l]{1}}
 \put(53.7,7.5){\makebox(2,1)[l]{4}}
 \put(52.5,2){\makebox(2,1)[l]{3}}

 \put(59.6,16.8){\makebox(2,1)[l]{3}}
 \put(60.7,11.5){\makebox(2,1)[l]{1}}
 \put(60.7,7.5){\makebox(2,1)[l]{4}}
 \put(59.5,2){\makebox(2,1)[l]{2}}

 \multiput(11,-4)(7,0){5}{\circle*{0.7}}  \put(10.6,-3.2){\makebox(2,1)[l]{4}}
 \multiput(9,-8)(7,0){5}{\circle*{0.7}}   \put(9,-10){\makebox(2,1)[l]{1}}
 \multiput(7,-12)(7,0){5}{\circle*{0.7}}  \put(6.5,-14){\makebox(2,1)[l]{2}}
 \multiput(13,-8)(7,0){5}{\circle*{0.7}}  \put(12.5,-10){\makebox(2,1)[l]{3}}
 \multiput(11,-4)(7,0){5}{\line(-1,-2){4}}
 \multiput(11,-4)(7,0){5}{\line(1,-2){2}}

  \put(17.6,-3.2){\makebox(2,1)[l]{4}}
  \put(16,-10){\makebox(2,1)[l]{1}}
  \put(13.5,-14){\makebox(2,1)[l]{3}}
  \put(19.5,-10){\makebox(2,1)[l]{2}}

  \put(24.6,-3.2){\makebox(2,1)[l]{4}}
  \put(23,-10){\makebox(2,1)[l]{2}}
  \put(20.5,-14){\makebox(2,1)[l]{3}}
  \put(26.5,-10){\makebox(2,1)[l]{1}}

  \put(31.6,-3.2){\makebox(2,1)[l]{2}}
  \put(30,-10){\makebox(2,1)[l]{4}}
  \put(27.5,-14){\makebox(2,1)[l]{1}}
  \put(33.5,-10){\makebox(2,1)[l]{3}}

  \put(38.6,-3.2){\makebox(2,1)[l]{2}}
  \put(37,-10){\makebox(2,1)[l]{4}}
  \put(34.5,-14){\makebox(2,1)[l]{3}}
  \put(40.5,-10){\makebox(2,1)[l]{1}}

 \multiput(46,-4)(7,0){2}{\circle*{0.7}}  \put(45.6,-3.2){\makebox(2,1)[l]{4}}
 \multiput(46,-8)(7,0){2}{\circle*{0.7}}  \put(46.7,-8.5){\makebox(2,1)[l]{2}}
 \multiput(44,-12)(7,0){2}{\circle*{0.7}} \put(43.5,-14){\makebox(2,1)[l]{1}}
 \multiput(48,-12)(7,0){2}{\circle*{0.7}} \put(47.5,-14){\makebox(2,1)[l]{3}}
 \multiput(46,-4)(7,0){2}{\line(0,-1){4}}
 \multiput(46,-8)(7,0){2}{\line(-1,-2){2}}
 \multiput(46,-8)(7,0){2}{\line(1,-2){2}}

  \put(52.6,-3.2){\makebox(2,1)[l]{1}}
  \put(53.7,-8.5){\makebox(2,1)[l]{4}}
  \put(50.5,-14){\makebox(2,1)[l]{2}}
  \put(54.5,-14){\makebox(2,1)[l]{3}}

 \end{picture}
 \vskip 3cm
 \caption{16 trees in $\mathcal{R}_{4,2}[\deg(4)>0]$}
 \label{fig-deg4}
 \end{figure}
Before we start our journey of constructing the bijection,
we present an inductive proof. In principle, it follows
from Theorem \ref{thm-main-comb} for the case $\deg_T(1)=1$.
For completeness, we include the inductive proof which is
slightly simpler than that of (\ref{tnk-r}).

{\em Inductive Proof of Theorem \ref{thm-bij-r}.}
For $n\geq 2$, the arguments of Shor \cite{Shor95} or Zeng
\cite{Zeng} imply the following identities:
\begin{eqnarray*}
\mbox{(i)}&& R_{n,k+1}[\deg(n)>0]=(n+k-1)R_{n-1,k}.
\qquad \qquad \qquad \qquad \qquad \qquad \qquad \\[5pt]
\mbox{(ii)}&&R_{n,k}[\deg(1)>0]\\[5pt] &&\quad =\,
  (n-2)R_{n-1,k}[\deg(1)>0]+R_{n-1,k}\\[5pt]
&&\quad \qquad  +\, (n+k-2)R_{n-1,k-1}[\deg(1)>0].\\[5pt]
\mbox{(iii)}&&R_{n,k}=(n-1)R_{n-1,k} + (n+k-2)R_{n-1,k-1}.
\end{eqnarray*}
Because of (i), $R_{n,k}[\deg(1)>0]=R_{n,k+1}[\deg(n)>0]$ can be
deduced from the following relation:
\begin{equation}
R_{n,k}[\deg(1)>0]=(n+k-1)R_{n-1,k}, \label{rnkk}
\end{equation}
for $n\geq 1$. The above claimed identity obviously holds for
$n=1$.
 Suppose (\ref{rnkk}) holds for $n-1$.
 From (i) -- (iii) and the inductive hypothesis, it follows that
\begin{eqnarray*}
\lefteqn{R_{n,k}[\deg(1)>0]}\\[5pt]
&=&(n-2)R_{n-1,k}[\deg(1)>0]+R_{n-1,k}+(n+k-2)R_{n-1,k-1}[\deg(1)>0]\\[5pt]
&=&(n-2)(n+k-2)R_{n-2,k}+R_{n-1,k}+(n+k-2)(n+k-3)R_{n-2,k-1}\\[5pt]
&=&(n+k-2)\bigl\{(n-2)R_{n-2,k}+(n+k-3)R_{n-2,k-1}
          \bigr \}+R_{n-1,k}\\[5pt]
&=&(n+k-2)R_{n-1,k}+R_{n-1,k}\\[5pt] &=&(n+k-1)R_{n-1,k}.
\end{eqnarray*}
This completes the proof. \qed

We note that for $k=n-1$, there does not exist any
rooted tree $T$ with $n$ nodes and $k$ improper edges such
that $\deg_T(1)>0$, because any edge with $1$ as the father
node is proper. Thus, we may assume without loss of generality
that $k<n-1$.

It turns out that we need to consider two major cases in the construction of
a bijection for (\ref{bij-r}).
First, we introduce the notation $\mathcal{R}_{n,k}^{(i)}$
for the set of trees $T$ in $\mathcal{R}_{n,k}$
such that there are $i$ proper edges on the path from the node $n$ to the
root. Suppose $T$ is a rooted tree on $[n]$ and $x$ is a node of $T$ such
that $T_x$ contains the node $n$. Then we may define the {\it lowering
operation} $L$ on $T_x$ such that $L(T_x)$ is the rooted tree obtained
from $T_x$ by taking $n$ as the new root and letting the ancestor
nodes of $n$ fall down to the descendants of $n$. Under certain
circumstances, the lowering operation is reversible, and the
reverse will be called the {\em lifting operation}. The following
Theorem \ref{thm-bij-r1} tells us where we may apply the lowering operations that are
reversible. We need to define the {\em upper critical node} and
the {\em lower critical node} of a rooted tree $T$.

If  $T$ is a rooted tree in $\mathcal{R}_{n,k}^{(i)}$, where
$i\geq 1$.  Suppose $(n=v_1, v_2, \ldots, v_t)$ is the path from
$n$ to the root of $T$, and $v_j$ is the first node on the path
such that $(v_{j-1}, v_j)$ is a proper edge of $T$. Then we call
$v_j$ the {\em upper critical node} of $T$.
 On the other hand, for any rooted tree $T$ on $[n]$
such that $\deg_T(n)>0$, we define the {\em lower critical node}
of $T$ by the following procedure. First, we note that
$\beta(n)<n$. By abuse of the above indices $t$ and $j$, we assume
that $(n=u_1, u_2, \ldots, u_t = \beta(n))$ is the path from $n$
to $\beta(n)$ and that $u_j\neq n$ is the first node on the path
such that every node in $T_n-T_{u_j}$, namely, the tree obtained
from $T_n$ by removing the subtree $T_{u_j}$, is greater than
$u_j$, denoted
\[ u_j < \beta(T_n -T_{u_j}). \]
Note that such a node $u_j$ must exist because the
node $\beta(n)$ is always a candidate to satisfy the above condition.
The lower critical node of $T$ will be denoted
by $\lambda(T)$, or $\lambda$ for short, if
no confusion arises in the context.

\setlength{\unitlength}{0.08in}
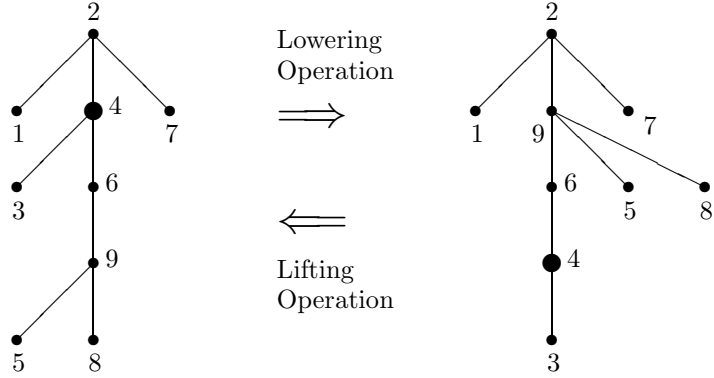
\begin{figure}[ht]
\hskip 1.3cm
\begin{picture}(60,20)
    \put(11,16){\circle*{0.7}}      \put(10.6,17){\makebox(2,1)[l]{2}}
    \put(11,11){\circle*{1.2}}      \put(12,10.8){\makebox(2,1)[l]{4}}
    \put(11,6){\circle*{0.7}}       \put(11.8,5.8){\makebox(2,1)[l]{6}}
    \put(11,1){\circle*{0.7}}       \put(11.8,0.8){\makebox(2,1)[l]{9}}
    \put(11,-4){\circle*{0.7}}      \put(10.7,-6){\makebox(2,1)[l]{8}}
    \put(6,-4){\circle*{0.7}}       \put(5.7,-6){\makebox(2,1)[l]{5}}
    \put(6,6){\circle*{0.7}}        \put(5.7,4){\makebox(2,1)[l]{3}}
    \put(6,11){\circle*{0.7}}       \put(5.7,9){\makebox(2,1)[l]{1}}
    \put(16,11){\circle*{0.7}}      \put(15.7,9){\makebox(2,1)[l]{7}}

  \put(11,16){\line(-1,-1){5}}
  \put(11,16){\line(0,-1){20}}
  \put(11,16){\line(1,-1){5}}
  \put(11,11){\line(-1,-1){5}}
  \put(11,1){\line(-1,-1){5}}

    \put(41,16){\circle*{0.7}}     \put(40.6,17){\makebox(2,1)[l]{2}}
    \put(41,11){\circle*{0.7}}     \put(39.7,9){\makebox(2,1)[l]{9}}
    \put(41,6){\circle*{0.7}}      \put(41.8,5.8){\makebox(2,1)[l]{6}}
    \put(41,1){\circle*{1.2}}      \put(42,0.8){\makebox(2,1)[l]{4}}
    \put(41,-4){\circle*{0.7}}     \put(40.7,-6){\makebox(2,1)[l]{3}}
    \put(36,11){\circle*{0.7}}     \put(35.7,9){\makebox(2,1)[l]{1}}
    \put(46,6){\circle*{0.7}}      \put(45.7,4){\makebox(2,1)[l]{5}}
    \put(46,11){\circle*{0.7}}     \put(47,9.5){\makebox(2,1)[l]{7}}
    \put(51,6){\circle*{0.7}}      \put(50.7,4){\makebox(2,1)[l]{8}}

  \put(41,16){\line(-1,-1){5}}
  \put(41,16){\line(0,-1){20}}
  \put(41,16){\line(1,-1){5}}
  \put(41,11){\line(1,-1){5}}
  \put(41,11){\line(2,-1){10}}

 \put(23,10){\makebox(2,1)[l]{\LARGE{$\Longrightarrow$}}}
 \put(23,13){\makebox(2,1)[l]{Operation}}
 \put(23,15){\makebox(2,1)[l]{Lowering}}
 \put(23,3){\makebox(2,1)[l]{\LARGE{$\Longleftarrow$}}}
 \put(23,0){\makebox(2,1)[l]{Lifting}}
 \put(23,-2){\makebox(2,1)[l]{Operation}}

 \end{picture}
 \vskip 1.5cm
 \caption{Lowering and Lifting Operations}
\end{figure}

With the aid of the lifting and lowering operations, we may establish the
following bijection which serves as the first case for the
bijection (\ref{bij-r}).

\begin{thm}\label{thm-bij-r1}
For $i\geq 1$, we have the following bijection:
\begin{equation}
     \mathcal{R}_{n,k}^{(i)}[\deg(1)>0]
     \longleftrightarrow
     \mathcal{R}_{n,k+1}^{(i-1)}[\deg(n)>0,
     (\deg(1)>0 \mbox{ or } \lambda=1)].
\label{bij-r1}
\end{equation}
\end{thm}

{\em Proof.} Suppose $T$ is tree in
$\mathcal{R}_{n,k}^{(i)}[\deg(1)>0]$. We assume that  $(n=v_1,
v_2, \ldots, v_t)$ is the path from $n$ to the root of $T$, and
$v_j$ is the upper critical node of $T$. We now apply the lowering
operation $L$ on $T({v_j})$. We then obtain a rooted tree $T'$ by
substituting the subtree $T({v_j})$ with $L(T({v_j}))$. Note that
the resulting tree $T'$ has one more improper edge than $T$
because the edge $(v_j, v_{j-1})$ is proper in $T$ and the edge
$(v_{j-1}, v_j)$ is improper in $T'$. Moreover, we notice
that after the lowering operation, the degree of $n$ increases by
$1$, the degree of the upper critical node decreases by 1, and the
degree of any other node remains unchanged. Therefore, if we have
$\deg_{T'}(1)=0$, then $1$ must be the upper critical node of $T$
because $\deg_T(1)>0$.

We now face the task of recovering the original tree
$T$ from the tree $T'$ and convincing ourselves of the fact that
the upper critical node of $T$ becomes the lower critical node of
$T'$.  In order to single out the upper critical node of $T$ in
the new environment of $T'$, we first claim that the upper
critical node of $T$, say $w$, has to be on the path from $n$ to
$\beta(n)$ in $T'$. Assume that $v$ is the child of $w$ that is on the
 path from $w$ to $n$. By the definition of $w$, one sees that
 $w<\beta(T_v)$. Therefore, after the application of
the lowering operation, $w$ has to be on the path from $n$ to
$\beta(n)$.

We now assume that $(n, u_1, u_2, \ldots)$
is the path from $n$ to $\beta(n)$ in $T'$.
If $u_1 < \beta(T'_n - T'_{u_1})$, then
one sees that $(u_1, n)$ is a proper edge in $T$ and
one can lift the edge $(u_1, n)$ up and to restore
$w$ as the upper critical node of $T$.
Otherwise, we may consider the next candidate $u_2$, and
so on. Such a process shows that the upper critical node
of $T$ can be identified by the lower critical node of $T'$.
This completes the proof.   \qed

The next case we should  consider is the following theorem.

\begin{thm}\label{thm-bij-r2}
For $n\geq 1$ and $m\geq 1$,
we have the following bijection:
\begin{equation}
     \mathcal{R}_{n,k}^{(0)}[\deg(1)=m]
    \longleftrightarrow
\mathcal{R}_{n,k+1}^{(m-1)}[\deg(n)>0,
     (\deg(1)=0\mbox{ and }\lambda>1 )].
\label{bij-r2}
\end{equation}
\end{thm}

Note that Theorems  \ref{thm-bij-r1} and \ref{thm-bij-r2} together
lead to a refined version of Theorem \ref{thm-bij-r}. We now focus
on the proof of (\ref{bij-r2}). The proof of (\ref{bij-r1})
actually implies the following assertion:

\begin{lem}\label{lem-rnk-fg1}
For $i\geq 1, m\geq 1$, we have
the following bijection:
\[  \mathcal{R}_{n,k}^{(i)}
[\deg(n)=m,\deg(1)=0, \lambda>1] \longleftrightarrow \qquad \qquad
\qquad \qquad
 \]
\begin{equation}
\qquad \qquad
  \mathcal{R}_{n,k+1}^{(i-1)}
[\deg(n)=m+1,\deg(1)=0,\lambda>1]. \label{rnk-fg1}
\end{equation}
\end{lem}

By iteration, for any $m\geq 1$ it follows that
\[ \mathcal{R}_{n,k+1}^{(m-1)}
[\deg(1)=0,\deg(n)\geq 1,\lambda>1] \qquad \qquad \qquad \qquad
\qquad \]
\begin{equation}
\qquad \qquad \qquad \longleftrightarrow \mathcal{R}_{n,k+m}^{(0)}
[\deg(1)=0,\deg(n)\geq m,\lambda>1]. \label{rnkm-1}
\end{equation}

Because of the above bijection, one sees that
Theorem \ref{thm-bij-r2} is equivalent to the
following statement.

\begin{thm}\label{thm-bij-rx}
For $n\geq 1$ and $m\geq 1$, we have the following
bijection:
\begin{equation}
      \mathcal{R}_{n,k}^{(0)}
[\deg(1)=m]\longleftrightarrow
 \mathcal{R}_{n,k+m}^{(0)}
[\deg(1)=0,\deg(n)\geq m,\lambda>1]. \label{bij-rx}
\end{equation}
\end{thm}

We now run short of notation and terminology for our
unaccomplished mission, and here are more in need.
\begin{itemize}
\item $\alpha\ ={\alpha}_{T}:=
\max\{{\beta}_{T}(b):\, b\mbox{ is a child
                 of the node }n \},$
   for $ T\in \mathcal{R}_{n,k}[\deg(n)>0].$
such that $\deg_{T}(n)>0.$
\item
${\beta}^{*}={\beta}^{*}_{_T}:=
\mbox{min}\{{\beta}_{_T}(a):\, a\mbox{ is a
                 child of the node } 1 \},$
for $T\in \mathcal{R}_{n,k}[\deg(1)>0].$

\item
$x\prec y$ denotes that $x$ is a descendant of $y$, while
$x\not\prec y$ means the opposite.

\item If we cut off a subtree from a node $u$ and
join it to another node $v$ as a subtree, we will
simply say that the subtree is moved to another node,
or we move the subtree to another node.
\end{itemize}

Note that, for any $T\in\mathcal{R}_{n,k}^{(0)}$, the node $1$
cannot be on the path from $n$ to the root, namely, $n\not\prec
1$. Also, if $T\in\mathcal{R}_{n,k}^{(0)}$, then we have
$\deg(n)>0$; Otherwise, the first edge on the path from $n$ to the
root would be proper. Therefore, $\alpha$ is always well-defined
for a tree $T\in\mathcal{R}_{n,k}^{(0)}$, and if $1\not\prec n$,
we have $\lambda_T>1$.

The following lemma is crucial.

\begin{lem}\label{lem-rnk-w}
We have the following bijection:
\begin{equation}
     \mathcal{R}_{n,k}[\deg(1)=1,\beta^*=w]\longleftrightarrow
     \mathcal{R}_{n,k+1}[\deg(1)=0,\mu=w],\
\label{rnk-w}
\end{equation}
where $\mu(T)$ is defined for any rooted tree in which  $1$ is not
the root of $T$.  We suppose that $(u_1=1, u_2, \ldots, u_t=v)$ is
the path from $1$ to the root of $T$. Then $\mu(T)=u_j$
denotes the first node on the path with $u_j>1$ such that
\begin{equation}
u_j< \beta(T_{v} -T_{u_j}).\label{cond-1}
\end{equation}
Moreover, we always assume that $v$ satisfies the above
condition (\ref{cond-1}).
\end{lem}

{\em Proof.} Suppose $T\in\mathcal{R}_{n,k}[\deg(1)=1,\beta^*=w]$,
where $w\geq 2$. Let $v$ be the
  unique child of $1$, and $P\colon \,(v_1,v_2,\ldots,v_t=1)$
the path from the root of $T$ to $1$. The scheme of the
construction consists of the following steps:
\begin{itemize}
\item Cut off the edge $(1, v)$ and get a tree
$S=T - T_{v}$.

\item Cut off some edges on the path from $v_1$ to $1$
to get a forest, say $R_1, R_2, \ldots, R_s$, subject
to some conditions to be spelled out later.

\item Obtain a tree $T'$ from $T_{v}$ by joining the
each $R_i$ as a subtree of the node $\beta_{T}(v)$
in $T_v$.
\end{itemize}

The tree $T'$ constructed above will be the goal of our
bijection. We now make it precise.

First, if $v_1< w$ then set $j_1=1$. Otherwise,
we choose $j_1$ as the minimum index such that
\begin{equation}
v_{j_1}< \beta(T(v_1)-T(v_{j_1})), \quad\mbox{ and }\quad
v_{j_1}<w. \label{cond-2}
\end{equation}
Because $v_t=1$ is on the path $P$, $j_1$
can be determined.
Second, we find all indices $j>j_1$
according to the following condition
\begin{equation}
v_{j}< \beta(T(v_1) -T(v_{j})), \label{cond-3}
\end{equation}
and denote by  $j_2,j_3,\ldots,j_s=t$, where  $j_2<j_3<\cdots<j_s$,
the solutions to the above inequality (\ref{cond-3}).
Third, set $j_0=0,\; v_{t+1}=v$ and
\[R_{i}=T(v_{j_{i-1}+1}) - T(v_{j_i+1}), \quad 1 \leq i\leq s, \]
namely,
\[R_1=T(v_1) -T(v_{j_1+1}), \quad
R_2= T(v_{j_1+1}) -T(v_{j_2+1}), \quad \ldots, \quad R_s=
T(v_{j_{s-1}+1}) -T(v).\]
Fourth, we construct a tree $T'$ from
the above decomposition of $T-T_v$ into $R_1, R_2, \ldots,R_s$ as
claimed before.

\setlength{\unitlength}{0.08in}
\begin{figure}[ht]
\hskip 2.7cm
\begin{picture}(60,22)

 \put(20,20){\circle*{0.7}} \put(21,19.5){\makebox(2,1)[l]{$v$}}
 \put(20,15){\circle*{0.7}}
 \put(20,10){\circle*{0.7}} \put(21,10){\makebox(2,1)[l]{$w$}}
 \put(20,20){\line(0,-1){15}}
 \put(20,10){\line(-2,-1){10}}
 \put(20,10){\line(-6,-5){6}}
 \put(20,10){\line(6,-5){6}}
 \put(20,10){\line(2,-1){10}}

  \multiput(9.9,2)(0,-1){3}{\makebox(2,1)[l]{$\cdot$}}
  \multiput(13.9,2)(0,-1){3}{\makebox(2,1)[l]{$\cdot$}}
  \put(15.5, 1.5){\makebox(2,1)[l]{$\ldots$}}
 \put(19,2){\makebox(2,1)[l]{$R_1 \;\ldots \;R_{s-1}\;R_s$}}

\end{picture}
\caption{$T'$}
\end{figure}

Let us now take a close look
at $T$ and $T'$.
Obviously, all the edges on the path from $v_1$ to
$v_t=1$ in $T$  are improper.
When we cut $T-T_v$ into $s$ pieces
$R_1, R_2, \ldots, R_s$, we lose
$s-1$ improper edges. By the definition of
$j_1,j_2,\ldots,j_s$, namely, the conditions (\ref{cond-2}) and
(\ref{cond-3}), one sees that $v_{j_1} < w$. Since
$j_2 > j_1$, the node $v_{j_1}$ must be a node
in $T(v_1)-T(v_{j_2})$. It follows that
\[ v_{j_2} < \beta(T(v_1) - T({v_{j_2}}))  \leq v_{j_1} .\]
The same reasoning leads to order relation:
\[ w>v_{j_1}>v_{j_2}>\cdots>v_{j_s}=1.\]
Therefore, when joining $R_1, R_2, \ldots, R_s$
as the subtrees of $w=\beta(T_v)$ in $T_v$,
we gain $s$ improper edges. Taking the previously lost improper edges
into consideration, one sees that $T'$ has one more improper edge
than $T$.

We now come to the justification of the fact that the node $w$ in
$T$ can be recovered as the node $\mu(T')$. From the above
construction, one sees that $w$ is a node on the path from the
node $1$ to the root. Moreover, $w$ satisfies condition
(\ref{cond-1}). We now need to show that except for $1$, there is
no other node $v_j$ on the path from $1$ to $w$ satisfying the
same condition (\ref{cond-1}). Suppose that $v_j$ is such a node,
namely, $v_j < \beta(T'(v) - T'(v_j))$. Since $w$ is the minimum
node in $T_v$, we have
\[ \beta(T'(v)-T'(v_j)) = \beta(T'(w)-T'(v_j))\]
Thus, we obtain
\begin{equation}
v_j< \beta(T(v_1)-T(v_j)), \quad\mbox{ and }\quad v_j<w.
\label{vjr}
\end{equation}
Let us consider what happens to the nodes on the path from $1$ to
the root in the above subtree $R_s=T(v_{j_{s-1}+1}) - T(v)$.
Suppose $v_j$ is such a node where $j_{s-1}+1 \leq j < j_s=t$. If
$s=1$, (\ref{vjr}) is contradictory to the definition of $j_1$. If
$s\geq 2$, since $j$ is not a solution to (\ref{cond-3}), we
have
\begin{equation}
v_j> \beta(T(v_1)-T(v_j)). \label{tvjsj}
\end{equation}
Note that two relations (\ref{vjr}) and (\ref{tvjsj})  are
contradictory to each other. Therefore, we reach the conclusion
that $T'\in\mathcal{R}_{n,k+1}[\deg(1)=0,\mu=w]$.

Before we give the reverse procedure to reconstruct $T$ from $T'$.
We need the following three claims about the above procedure.

\noindent
{\em Claim 1.} The node $v_{j_i}$ must be
on the path from $v_1$ to $\beta(R_i)$ in $R_i$.

The condition (\ref{cond-2}) says that every node in
$R_1'=T(v_1)-T(v_{j_1})$ is greater than $v_{j_1}$. The subtree
$R_1$ consists of $R_1'$ joined by a subtree rooted at $v_{j_1}$.
Thus, $v_{j_1}$ must be on the path from $v_1$ to $\beta(R_1)$.

For $R_2, R_3, \ldots, R_s$, the same argument applies.
Since $j_{i-1}+1\leq j_i,\; j_i+1>j_i$. It follows that
$v_{j_i}\in T(v_{j_{i-1}+1})$ and $v_{j_i}\not\in T(v_{j_i+1})$,
that is $v_{j_i}\in R_i=T(v_{j_{i-1}+1}) - T(v_{j_i+1})$. By the
definition of $v_{j_{i+1}}$, we have
\[v_{j_i} < \beta(T(v_1)-T(v_{j_i})) \leq \beta(T(v_{j_{i-1}+1})-T(v_{j_i}))
 = \beta(R_i-R_i(v_{j_i})).\]
Thus, $\beta(R_i)$ is in $R_i(v_{j_i})$. Note that we have assumed that the
above equalities is true for $v_{j_{i-1}+1}=v_{j_i}$.

\noindent
{\em Claim 2.}
$w> \beta(R_1) > \beta(R_2) > \cdots > \beta(R_s)$.

From the following relations
\[v_{j_i} < \beta(T(v_1)-T(v_{j_i})) \leq
\beta(T(v_{j_{i-2}+1})-T(v_{j_{i-1}+1})) = \beta(R_{i-1}), \]
we obtain $\beta(R_i) \leq v_{j_i} < \beta(R_{i-1})$.
We have already shown that
 $\beta(R_1) \leq v_{j_1} < w $, so Claim 2 follows.

\noindent
{\em Claim 3.}
The node $v_{j_i}$ can be determined
as the first node $z_j$ on the path
$z_1,z_2,\ldots,z_r=\beta(R_i)$ from the root of $R_i$ to
$\beta(R_i)$ such that
\begin{equation}
 z_j < \beta (R_i(z_1) - R_i(z_{j}))\quad \mbox{ and }\quad
 z_j < \beta(R_{i-1}).\label{cond-4}
\end{equation}
Here, we set $\beta(R_0)=w$ and assume the first inequality of
(\ref{cond-4}) is always true for $z_j=z_1$.

It is easy to see  that Claim 3 holds for $i=1$.
We now assume that $i\geq 2$.
In the proofs of Claim 1 and Claim 2, we have also shown that
\[v_{j_i} < \beta(R_i-R_i(v_{j_i}))\quad \mbox{ and }
\quad v_{j_i} < \beta(R_{i-1}).\]
Suppose there is another $v_j\neq v_{j_i}$ on the path from the
root of $R_i$ to $v_{j_i}$ that satisfies the condition
\[v_j < \beta(R_i-R_i(v_j))\quad \mbox{ and }
\quad v_j < \beta(R_{i-1}).\]
 By Claim 2 and the above condition, $v_j$ satisfies the
condition
\[ v_{j}< \beta(T(v_1) -T(v_{j})). \]
 However, we must have $j_{i-1}+1 \leq j <j_i$,
because $v_{j_{i-1}+1}$ is the root of $R_i$.
This is a contradiction to the fact that $j_2<\cdots<j_s$
are the only solutions greater than $j_1 $ to the above inequality.

We now come to the turning point of the bijection. For a tree
$T'\in \mathcal{R}_{n,k+1}[\deg(1)\linebreak
=0,\mu=w]$, we are
going to reconstruct the tree $T$. The first step is easy: the
subtrees $R_1, R_2, \ldots, R_s$ can be separated from $T'$ as the
subtrees $R_i$ of $w$ such that $\beta(R_i) < w$. By Claim 1,
$R_1, R_2, \ldots, R_s$ can be restored by the following order:
\[w> \beta(R_1) > \beta(R_2) > \cdots  > \beta(R_s).\]

Let $R=T'-R_1-R_2-\cdots -R_s$.
By the construction of $T'$, we need to merge the
subtrees $R_1, R_2, \ldots, R_s$ into a rooted
tree $S$. In so doing, we need to identify
which node on $R_i$ is the last node on the path
from $v_1$ to the node $1$ in $T$. In other words,
we need to have the nodes $v_{j_1}$,
$v_{j_2}$, $\ldots$, $v_{j_s}$ restored. This
job can be left to Claim 2.

The last step would be to put $R_1$, $R_2$, $\ldots$,
$R_s$ together with $R$.
For $i=1, 2,\ldots,  s-1$,
we join $R_{i+1}$ as a subtree of $v_{j_i}$ of
$R_i$, then join $R$ as a subtree of $1$ in $R_s$.
At last, we obtain the tree $T\in\mathcal{R}_{n,k}[\deg(1)=1,\beta^*=w]$.
\qed

Here is an example for $n=20,\, w=11$.

\setlength{\unitlength}{0.07in}
\begin{figure}[ht]
\hskip 2.2cm
\footnotesize{
\begin{picture}(60,20)
    \put(-4,14){\circle*{0.7}}       \put(-4.9,14.8){\makebox(2,1)[l]{9}}
    \put(-6,10){\circle*{0.7}}       \put(-5.3,9.8){\makebox(2,1)[l]{15}}
    \put(-2,10){\circle*{0.7}}       \put(-1.3,9.8){\makebox(2,1)[l]{16}}
    \put(-8,6){\circle*{0.7}}        \put(-7.3,5.8){\makebox(2,1)[l]{10}}
    \put(-4,6){\circle*{0.7}}        \put(-3.3,5.8){\makebox(2,1)[l]{12}}
    \put(0,6){\circle*{0.7}}         \put(0.7,5.8){\makebox(2,1)[l]{7}}
    \put(-10,2){\circle*{0.7}}       \put(-10.7,0){\makebox(2,1)[l]{20}}
    \put(-6,2){\circle*{0.7}}        \put(-6.7,0){\makebox(2,1)[l]{18}}
    \put(-2,2){\circle*{0.7}}        \put(-2.7,0){\makebox(2,1)[l]{6}}
    \put(2,2){\circle*{0.7}}         \put(2.7,1.8){\makebox(2,1)[l]{8}}
    \put(4,-2){\circle*{0.7}}        \put(4.7,-2.2){\makebox(2,1)[l]{3}}
    \put(2,-6){\circle*{0.7}}        \put(1.7,-8){\makebox(2,1)[l]{4}}
    \put(6,-6){\circle*{0.7}}        \put(6.7,-6.2){\makebox(2,1)[l]{17}}
    \put(4,-10){\circle*{0.7}}       \put(3.7,-12){\makebox(2,1)[l]{2}}
    \put(8,-10){\circle*{0.7}}       \put(8.7,-10.2){\makebox(2,1)[l]{5}}
    \put(10,-14){\circle*{0.7}}      \put(10.7,-14.2){\makebox(2,1)[l]{1}}
    \put(12,-18){\circle*{0.7}}      \put(12.7,-18.2){\makebox(2,1)[l]{14}}
    \put(10,-22){\circle*{1.2}}      \put(9.9,-24){\makebox(2,1)[l]{11}}
    \put(14,-22){\circle*{0.7}}      \put(13.7,-24){\makebox(2,1)[l]{13}}
    \put(8,-26){\circle*{0.7}}       \put(7.1,-28){\makebox(2,1)[l]{19}}

  \put(-4,14){\line(-1,-2){6}}
  \put(-4,14){\line(1,-2){18}}
  \put(-2,10){\line(-1,-2){4}}
  \put(0,6){\line(-1,-2){2}}
  \put(4,-2){\line(-1,-2){2}}
  \put(6,-6){\line(-1,-2){2}}
  \put(12,-18){\line(-1,-2){4}}

    \put(-4,-18){\makebox(2,1)[l]{$T$}}
    \put(41,-18){\makebox(2,1)[l]{$T'$}}\hskip -0.6cm
    \put(41,6){\circle*{0.7}}       \put(40.1,6.8){\makebox(2,1)[l]{14}}
    \put(38,2){\circle*{1.2}}       \put(39.3,1.9){\makebox(2,1)[l]{11}}
    \put(49,2){\circle*{0.7}}       \put(48.6,0){\makebox(2,1)[l]{13}}
    \put(30,-2){\circle*{0.7}}       \put(29.1,-4){\makebox(2,1)[l]{19}}
    \put(35,-2){\circle*{0.7}}       \put(34.7,-4){\makebox(2,1)[l]{9}}
    \put(38,-2){\circle*{0.7}}       \put(38.7,-2.2){\makebox(2,1)[l]{16}}
    \put(42,-2){\circle*{0.7}}       \put(41.6,-4){\makebox(2,1)[l]{8}}
    \put(46,-2){\circle*{0.7}}       \put(45.6,-1.2){\makebox(2,1)[l]{17}}
    \put(32,-6){\circle*{0.7}}       \put(31.7,-8){\makebox(2,1)[l]{15}}
    \put(36,-6){\circle*{0.7}}       \put(36,-8){\makebox(2,1)[l]{12}}
    \put(40,-6){\circle*{0.7}}       \put(40.7,-6.2){\makebox(2,1)[l]{7}}
    \put(46,-6){\circle*{0.7}}       \put(45.7,-8){\makebox(2,1)[l]{3}}
    \put(50,-6){\circle*{0.7}}       \put(49.7,-8){\makebox(2,1)[l]{2}}
    \put(54,-6){\circle*{0.7}}       \put(53.7,-5.2){\makebox(2,1)[l]{5}}
    \put(29,-10){\circle*{0.7}}      \put(28.7,-12){\makebox(2,1)[l]{10}}
    \put(34,-10){\circle*{0.7}}       \put(33.1,-12){\makebox(2,1)[l]{18}}
    \put(42,-10){\circle*{0.7}}       \put(41.7,-12){\makebox(2,1)[l]{6}}
    \put(62,-10){\circle*{0.7}}       \put(61.7,-12){\makebox(2,1)[l]{1}}
    \put(50,-10){\circle*{0.7}}       \put(49.7,-12){\makebox(2,1)[l]{4}}
    \put(26,-14){\circle*{0.7}}       \put(25.1,-16){\makebox(2,1)[l]{20}}

  \put(41,6){\line(-3,-4){3}}
  \put(41,6){\line(2,-1){8}}
  \put(38,2){\line(-2,-1){8}}
  \put(38,2){\line(0,-1){4}}
  \put(38,2){\line(-3,-4){12}}
  \put(38,2){\line(1,-1){12}}
  \put(38,2){\line(2,-1){24}}
  \put(38,-2){\line(-1,-2){4}}
  \put(38,-2){\line(1,-2){4}}
  \put(46,-2){\line(1,-1){4}}
 \put(16,-2){\makebox(2,1)[l]{\LARGE{$\Longleftrightarrow$}}}
\end{picture}
}
\vskip 5cm \caption{Example for $n=20$ and $w=11$}
\end{figure}
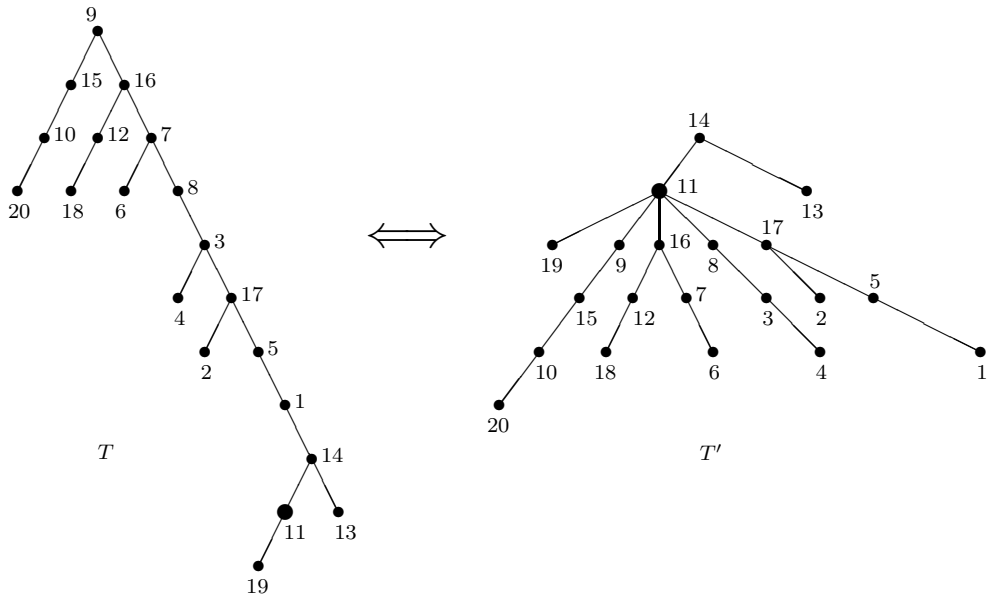

 We are now ready to present the proof of Theorem
\ref{thm-bij-rx} in the following refined version.
\begin{thm}\label{abcde}
For $m\geq 1$, we have the following
bijections:
\begin{description}
\item[{\bf (a)}]
$\mathcal{R}_{n,k}^{(0)}[\deg(1)=m,1\not\prec n,\alpha<\beta^{*}]$

\hskip 3cm $\longleftrightarrow \mathcal{R}_{n,k+m}^{(0)}
         [\deg(1)=0,\deg(n)\geq m+1,1\not\prec n]$.

\item[{\bf (b)}] $ \mathcal{R}_{n,k}^{(0)}
[\deg(1)=m,\alpha>\beta^{*}]$

   \hskip 3cm $\longleftrightarrow
          \mathcal{R}_{n,k+m}^{(0)} [\deg(1)=0,\deg(n)=m,1\not\prec n]$.

\item[{\bf (c)}]
$\mathcal{R}_{n,k}^{(0)}[\deg(1)=m,1\prec n,
         \deg(n)\geq 2,\alpha<\beta^{*}]$.

\hskip 3cm $\longleftrightarrow\mathcal{R}_{n,k+m}^{(0)}
         [\deg(1)=0,\deg(n)\geq m+1,1\prec n,\lambda>1]$.

\item[{\bf (d)}]
$\mathcal{R}_{n,k}^{(0)}
[\deg(1)=m,1\prec n,\deg(n)=1]$

\hskip 3cm $\longleftrightarrow\mathcal{R}_{n,k+m}^{(0)}
   [\deg(1)=0,\deg(n)=m,1\prec n,\lambda>1]$.
\end{description}
\end{thm}

{\em Proof.} Recall the facts that for any tree
$T\in\mathcal{R}_{n,k}^{(0)}$ we have $\deg(n)>0$ and that
$\alpha(T)$ is well defined. Moreover, we do not get into the
detailed discussion about the range of $m$, because when $m$ is
out of range the bijection would simply do nothing.

{\bf (a)}
Suppose $T\in\mathcal{R}_{n,k}^{(0)}
                  [\deg(1)=m,1\not\prec n,\alpha<\beta^{*}]$.

Since $n$ is not on the path from $1$ to the root and $1$ is not
on the path from $n$ to the root, $1$ and $n$ lie in different
branches of their minimum common ancestor, in other words, the
common ancestor furthest from the root.
 Moving all subtrees of
$1$ to the node $n$, we are led to a tree\\
\centerline{
$T'\in\mathcal{R}_{n,k+m}^{(0)}
     [\deg(1)=0,\deg(n)\geq m+1,1\not\prec n]$.}

Conversely, given the tree $T'$,  we assume that
     $b_1,b_2,\ldots,b_j(j\geq m+1)$ are
 the children of $n$ ordered by
$\beta(b_1)>\beta(b_2)>\cdots>\beta(b_j)$. We now move the first
$m$ subtrees $T'_{b_i}(1\leq i\leq m)$  to node $1$.
     Thus, we have recovered the above tree $T$.

 {\bf (b)} Suppose $T\in\mathcal{R}_{n,k}^{(0)}[\deg(1)=m,\alpha>\beta^{*}]$.
 Exchange the node $n$ and the subtree $T_{\alpha}$. Thus  the
 degree of $n$ becomes zero, the edges on the path from the root to
 $\alpha$ are all improper by the definition of $\mathcal{R}_{n,k}^{(0)}$,
 while the first edge on the path from $\alpha$ to $n$ is proper by the
 definition of $\alpha_T$.
 Then move all subtrees of $1$ to the node $n$.
Since $\alpha>\beta^{*}$,
 we obtain a tree\\
\centerline{ $T'\in\mathcal{R}_{n,k+m}^{(0)}
              [\deg(1)=0,\deg(n)=m,1\not\prec n]$.}

The reverse of the above procedure is strictly the other way around.
Starting with the above tree $T'$,
 move all the subtrees of $n$ back to the node $1$. Suppose
that we obtain a tree $T''$ and that
the path from the root to $n$ in $T''$ is
 $P:(y_1,y_2,\ldots,y_s=n)$. Let $(y_i,y_{i+1})$
 be the first proper edge on the path $P$. Suppose $R_1,R_2,\ldots,R_s$
 are all of the subtrees of $y_i$ such that $\beta (R_j)>y_i$, and
 $n\not\in R_j,\,\forall j=1,2,\ldots,s$. Move these subtrees
to the node  $n$, and exchange labels of
the nodes $y_i$ and $n$. Therefore, we get the above tree $T$
such that $y_i=\alpha (T)$.

 {\bf (c)} Suppose $T\in\mathcal{R}_{n,k}^{(0)}
     [\deg(1)=m,1\prec n,\deg(n)\geq 2,\alpha<\beta^{*}]$.
      Assume that  $b_1,b_2,\ldots,b_j(j\geq 2)$
are the children of $n$ ordered by
$ 1=\beta(b_1)<\beta(b_2)<\cdots<\beta(b_j)<{\beta}^*.$
Suppose $Q:\;(b_2=c_1,c_2,\ldots,c_t=\beta(b_2))$
is the path from $b_2$ to $\beta(b_2)$.
We locate  the first $c_i$ such that
\[ c_i<\beta(T(b_2)-T(c_i)),\quad c_i<{\beta}^*, \quad \mbox{and}
\quad c_i<\beta(b_3), \quad \mbox{if $j\geq 3$.}\] Moving the
subtree $T_{b_1}$ to the node $c_i$ and moving all subtrees of $1$
to the node $n$, we obtain a tree\\ \centerline{
$T'\in\mathcal{R}_{n,k+m}^{(0)}
             [\deg(1)=0,\deg(n)\geq m+1,1\prec n,\lambda>1]$,}
with the property that $c_i=\lambda(T')$.

 Conversely, for this tree $T'$, we have $x=\lambda(T')>1$. We assume that
     $d_1,d_2,\ldots,d_s(s\geq m+1)$ are the  children of $n$
ordered by $ \beta(d_1)>\beta(d_2)>\cdots>\beta(d_s)=1.$ Moving
$T'_{d_i}(1\leq i\leq m)$ to the node $1$, and moving the subtree
of $x$ that contains $1$ to the node $n$, we get the above tree $T$.

{\bf (d)} Suppose $T\in\mathcal{R}_{n,k}^{(0)}
      [\deg(1)=m,1\prec n,\deg(n)=1]$, and let $b$ be the unique
      child of $n$.
Assume that $a_1,a_2,\ldots,a_m$ are the children
of $1$ ordered by $\beta(a_1)<\beta(a_2)<\cdots<\beta(a_m).$
Moving $T_{a_i} \;(2\leq i\leq m)$
to the node $n$, and let $R$ be the resulting tree.
 Substituting $S=R_b$
      with $S'$ by applying  Lemma \ref{lem-rnk-w}, we obtain a tree\\
 \centerline{$T'\in\mathcal{R}_{n,k+m}^{(0)}
      [\deg(1)=0,\deg(n)=m,1\prec n,\lambda>1]$.}

 Conversely, for this tree $T'$, we assume that
the subtree of $n$ that contains the node $1$ is $S'$. Applying
Lemma \ref{lem-rnk-w} to $S'$, we may recover $S$. Moving the
other $m-1$ subtrees of $n$ to the node $1$,
we obtain the above tree $T$.            \qed

After such an exciting and exhausting journey, we finally come to
our destination---Theorem \ref{thm-bij-r}. The essence of the
Theorem \ref{thm-bij-r} is the duality between the minimum element
and the maximum element in a rooted tree. It is easy to imagine
that the labels of a rooted tree do not have to be a consecutive
segment of integers in order for the bijection to hold. For this
reason, we say that a rooted tree $T$ is relabeled by a set $V$ of
the same number of nodes if the minimum node of $T$ is relabeled
by the minimum node in $V$, the second minimum node is relabeled
by the second minimum node in $V$, and so forth. By applying
Theorem \ref{thm-bij-r}, we can construct the main bijection of
this paper,  leading to a combinatorial proof of Theorem
\ref{thm-main-comb}.

\begin{thm}\label{main-b-r}
For $1\leq r \leq n$ and  $0\leq k<n-r$, we have the following bijection:
\begin{equation}
   \mathcal{T}_{n+1,k}[\deg(2)>0,\deg(1)=r]
   \longleftrightarrow\mathcal{T}_{n+1,k+1}
[\deg(n+1)>0,\deg(1)=r].
\end{equation}
\end{thm}

{\em Proof.}  It is obvious that for $r=n$ both sides of
(\ref{main-b-r}) are empty. So we may assume that $1\leq r \leq
n-1$. First, it is easy to see that the case $r=1$ reduces to
Theorem \ref{thm-bij-r}, by relabeling the set $\{2, 3, \ldots,
n+1\}$ with $\{1, 2, \ldots, n\}$. Thus, we may assume that $r\geq
2$. Suppose $T\in\mathcal{T}_{n+1,k} [\deg(2)>0,\deg(1)=r]$.
Assume $x$ is the child of the root $1$ such that $2$ is a descendant
of $x$ in $T$, and $y$ is the child of the root $1$ such that
$n+1$ is a descendant of $y$ in $T$. Note that it is possible that
$x=y$. We now proceed to construct a tree $T' \in
\mathcal{T}_{n+1,k+1}
     [\deg(n+1)>0,\deg(1)=r]$.
We have three cases to consider:

{\em Case 1.}  $x=y$. In this case, the minimum element $2$ and
the maximum element $n+1$ both appear in the subtree $T_x$.
Applying the Theorem \ref{thm-bij-r} on $T_x$, we are led to a
rooted tree $T_x'$. Substituting $T_x$ by $T_x'$ in $T$, we obtain
a rooted tree $T' \in  \mathcal{T}_{n+1,k+1}
     [\deg(n+1)>0,\deg(1)=r]$.

{\em Case 2.} $x\neq y$, and $\deg_{T}(n+1)>0$.
We also apply Theorem \ref{thm-bij-r}
on $T_x$, and we may obtain a rooted tree
$T' \in  \mathcal{T}_{n+1,k+1}
     [\deg(n+1)>0,\deg(1)=r]$ as in Case 1.

 {\em Case 3.} $x\neq y$, and $\deg_{T}(n+1)=0$. Let us relabel
the subtrees $T_x$ and $T_y$. Suppose $T_x$ has nodes $2$ and
$u_1< u_2< \cdots < u_i$ and $T_y$ has nodes $n+1$ and $v_1< v_2<
\cdots < v_j$.
 Let $R$ be the rooted tree obtained
from $T_x$ by relabeled by $u_1< u_2 < \cdots < u_i$ and $n+1$,
and $S$ be the rooted tree obtained from $T_y$ relabled by $2$ and
$v_1 < v_2 < \cdots < v_j$. Applying Theorem \ref{thm-bij-r} on
$R$, we obtain a rooted tree $R'$ with $\deg_{R'}(n+1)>0$. Now
substituting $T_x$ by $R'$ and $T_y$ by $S$, we are led to a
rooted tree $T'$, which is clearly in $\mathcal{T}_{n+1,k+1}
     [\deg(n+1)>0,\deg(1)=r]$.

Since all the above steps are reversible. We now only need to
classify the cases for a tree $T' \in  \mathcal{T}_{n+1,k+1}
     [\deg(n+1)>0,\deg(1)=r]$ so that
they can fit into one of the above three cases.

{\em Case A.} If $2$ and $n+1$ are in the same subtree $T_x'$
where $x$ is a child of the root $1$, then we
resort to the reverse of Case 1 to recover the
tree $T\in\mathcal{T}_{n+1,k}
[\deg(2)>0,\deg(1)=r]$.

{\em Case B.} Suppose $u$ and $v$ are the children of $1$ in $T'$
such that $T_u'$ contains $n+1$ and $T_v'$ contains $2$. If the
degree of the maximum element in $T_v'$ is nonzero, then we may
resort the reverse of the construction in Case 2 to recover $T$.
Otherwise, the degree of the maximum element in $T_v'$ equals
zero. In this case, we may count on the reverse procedure of Case
3 to recover the desired $T$.  This completes the proof. \qed

\newsection{Open Problems}

In evaluation of the bijections presented in this paper, the
construction of (\ref{bij-r2}) seems to be much more technical
than it should be, especially when compared with the case
(\ref{bij-r1}). We would very  much like to propose the following
problem.
\begin{problem}
Find an intrinsic construction for the bijection
(\ref{bij-r}).
In particular, Lemma \ref{lem-rnk-w}
deserves a better explanation.
\end{problem}

The inductive proofs
of Theorem \ref{thm-main-comb} and Theorem \ref{thm-bij-r}
might serve a hint  if they can be informatively translated into
a bijective scheme.

The next problem is concerned with a refined
version of the recurrence relation for the
numbers $f_{n,k}=|\mathcal{R}_{n,k}|$.
Recall that $\lambda(T)$ denotes the lower critical node of $T$,
as defined in the previous section. Notice
that $\lambda(T)$ is defined only for a tree $T$ such
that $\deg_T(n)>0$, where $n$ is the maximum node.
For notational simplicity, we leave out the
condition $\deg_T(n)>0$ when the condition
$\lambda=i$ is present.
We have the following conjecture:

\begin{conj}\label{conj-fnk}
For $n\geq 3$ and $1\leq i\leq n-2$, we have the recurrence relation:
\begin{equation}
|\mathcal{R}_{n,k}[\lambda=i] |=
(n-2)\,|\mathcal{R}_{n-1,k}[\lambda=i] | +(n+k-3)
\, |\mathcal{R}_{n-1,k-1}[\lambda=i] |.
\label{fnk}
\end{equation}
\end{conj}

Some numerical evidence in support of the above conjecture
is presented below for speculation.

 \[\mbox{  Table of $R_{n,k}[\lambda=1]$.\quad\quad\quad
           Table of $R_{n,k}[\lambda=2]$.\quad\quad\quad
           Table of $R_{n,k}[\lambda=3]$.  }\]

 \centerline{  \begin{tabular}{|l|c|c|c|c|}\hline
    $ k\backslash n$  & $2$ & $3$ & $4$  & $5$  \\ \hline
      $1$             & $1$ & $1$ & $2$  & $6$  \\ \hline
      $2$             &     & $2$ & $7$  & $29$ \\ \hline
      $3$             &     &     & $8$  & $59$ \\ \hline
      $4$             &     &     &      & $48$ \\ \hline
   \end{tabular}
   \quad\quad
   \begin{tabular}{|l|c|c|c|c|}\hline
    $ k\backslash n$  & $2$ & $3$ & $4$  & $5$  \\ \hline
      $1$             &     & $1$ & $2$  & $6$  \\ \hline
      $2$             &     & $1$ & $5$  & $23$ \\ \hline
      $3$             &     &     & $4$  & $37$ \\ \hline
      $4$             &     &     &      & $24$ \\ \hline
   \end{tabular}
   \quad\quad
   \begin{tabular}{|l|c|c|c|c|}\hline
    $ k\backslash n$  & $2$ & $3$ & $4$  & $5$  \\ \hline
      $1$             &     &     & $2$  & $6$  \\ \hline
      $2$             &     &     & $4$  & $20$ \\ \hline
      $3$             &     &     & $3$  & $29$ \\ \hline
      $4$             &     &     &      & $18$ \\ \hline
   \end{tabular}  }
\vskip 0.5cm

Here are some very special cases:
\begin{eqnarray}
     &
|\mathcal{R}_{n,1}[\lambda=i] |
= |\mathcal{R}_{n-1,0}|=(n-2)!,  \ \ 1\leq i\leq n-1,
\label{s-1}\\
\equsep
     &
|\mathcal{R}_{n,k}[\lambda=n-1] |=
|\mathcal{R}_{n-1,k-1}|, \ \ 1\leq k\leq n-1.
\label{s-2}
\end{eqnarray}

If the above conjecture is true,
then we can use induction to derive
the recurrence relation (\ref{shor-rec-n})
     from (\ref{fnk}), (\ref{s-2}), and the obvious identity
\[ |\mathcal{R}_{n,k}[\deg(n)=0] | =(n-1)
     |\mathcal{R}_{n-1,k}|. \]

The following special case is worth mentioning:
\begin{equation}
     |\mathcal{R}_{n,n-1}|=(2n-3)!! = 1\cdot 3 \cdot \cdots \cdot (2n-3).
\end{equation}

We may make a connection to increasing plane trees. A rooted
tree on $[n]$ is called {\em increasing} if any path from
the root to another vertex forms an increasing sequence.
As  an equivalent statement, we may say that an increasing
tree is a rooted tree without improper edges. We have
the following observation

\begin{prop}\label{prop-p}
There is a bijection between the set of rooted trees on $[n]$ with
$n-1$ improper edges and the set of plane trees on $[n]$ without
improper edges.
\end{prop}

Note that $(2n-3)!!$ is also the number of increasing plane trees on
$[n]$. Here is a combinatorial interpretation. Let $T$ be a tree
in $\mathcal{R}_{n,n-1}$. Then $1$ has to be a leaf of $T$, and
all edges of $T$ are improper. Suppose that $(1, v_1, v_2, \ldots,
v_t)$ is the path from $1$ to the root. Then we may
recursively construct an increasing plane tree $T'$. If $T$ has
only one node, then $T'$ is the same as $T$. If $n>1$, then for
each $T(v_i)$ construct the corresponding increasing plane tree,
and put them together by joining them to the minimum node
in the order of $(v_1, v_2, \ldots, v_t)$. An example is given
in Figure 7.  The tree on the left is a
rooted tree in which every edge is improper, and
the tree on the right is an increasing plane tree.

\setlength{\unitlength}{0.08in}
\begin{figure}[ht]
\hskip 0.5cm
\begin{picture}(60,20)

    \put(16,16){\circle*{0.7}}      \put(15.6,16.8){\makebox(2,1)[l]{4}}
    \put(12,12){\circle*{0.7}}      \put(11.3,12.8){\makebox(2,1)[l]{6}}
    \put(8,8){\circle*{0.7}}        \put(7.3,8.8){\makebox(2,1)[l]{9}}
    \put(4,4){\circle*{0.7}}        \put(3.5,2){\makebox(2,1)[l]{1}}
    \put(8,4){\circle*{0.7}}        \put(7.5,2){\makebox(2,1)[l]{5}}
    \put(12,4){\circle*{0.7}}       \put(11.7,2){\makebox(2,1)[l]{8}}
    \put(16,8){\circle*{0.7}}       \put(15.5,6){\makebox(2,1)[l]{2}}
    \put(20,12){\circle*{0.7}}      \put(19.4,10){\makebox(2,1)[l]{7}}
    \put(24,8){\circle*{0.7}}       \put(23.5,6){\makebox(2,1)[l]{3}}

  \put(16,16){\line(-1,-1){12}}
  \put(16,16){\line(1,-1){8}}
  \put(12,12){\line(1,-1){4}}
  \put(8,8){\line(0,-1){4}}
  \put(8,8){\line(1,-1){4}}

    \put(51,16){\circle*{0.7}}     \put(50.6,16.8){\makebox(2,1)[l]{1}}
    \put(47,12){\circle*{0.7}}     \put(46.3,12.8){\makebox(2,1)[l]{5}}
    \put(43,8){\circle*{0.7}}      \put(42.3,8.8){\makebox(2,1)[l]{8}}
    \put(39,4){\circle*{0.7}}      \put(38.5,2){\makebox(2,1)[l]{9}}
    \put(51,12){\circle*{0.7}}     \put(51.7,11){\makebox(2,1)[l]{2}}
    \put(51,8){\circle*{0.7}}      \put(50.5,6){\makebox(2,1)[l]{6}}
    \put(55,12){\circle*{0.7}}     \put(55,12.8){\makebox(2,1)[l]{3}}
    \put(55,8){\circle*{0.7}}      \put(54.5,6){\makebox(2,1)[l]{7}}
    \put(59,8){\circle*{0.7}}      \put(58.5,6){\makebox(2,1)[l]{4}}

  \put(51,16){\line(-1,-1){12}}
  \put(51,16){\line(0,-1){8}}
  \put(51,16){\line(1,-1){8}}
  \put(55,12){\line(0,-1){4}}

\put(29,9){\makebox(2,1)[l]{\LARGE{$\Longleftrightarrow$}}}
\end{picture}
\caption{Example for Proposition \ref{prop-p}}
\end{figure}
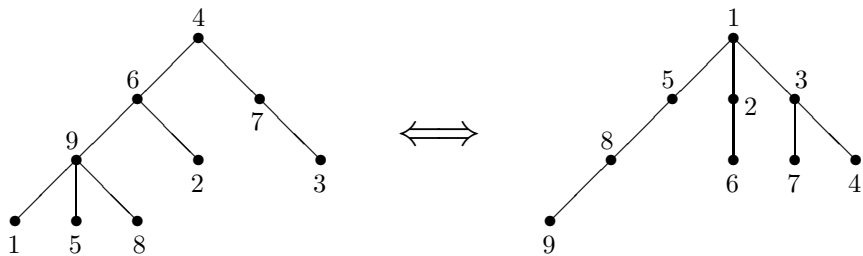

We now state a problem based on the above simple
observation, yet to be better understood.

\begin{problem}
Since $Q_{n,0}(x)$ corresponds to increasing trees on
$[n]$ while $Q_{n,n-1}(x)$ corresponds to increasing plane trees on
$[n]$, there must be some
kind of combinatorial structure like partial increasing plane
trees, a notion of interpolation of increasing trees and
increasing plane trees. Such a structure should
serve the purpose as an alternative combinatorial
interpretation of the Ramanujan polynomials.
\end{problem}

It is quite intriguing that there lie rich combinatorial
structures behind the Ramanujan polynomials. No doubt
that we may expect more episodes of uncovering further
mysteries plotted by these polynomials. Hopefully, we have made
some room for imagination, and we may (in any case) keep our fingers crossed
with respect to further developments.

\vskip 0.3cm
\par {\bf Acknowledgments.}
We would like
to thank J. Zeng for introducing us to the fascinating
subject of Ramanujan polynomials during his stimulating lectures given
at Nankai University in 1996.  We would also like to thank D. C. Torney
and C. Wang for helpful comments and discussions.


\begin{thebibliography}{99}
\small \setlength{\itemsep}{-.8mm}

\bibitem{Berndt83}  B. C. Berndt,  R. J. Evans and B. M. Wilson,  Chapter 3 
of Ramanujan's second notebooks, Adv. Math. {\bf 49} (1983), 123-169.

\bibitem{Berndt85}
B. C. Berndt,  Ramanujan's Second Notebooks, Part   I, Chap. 3: 
Combinatorial analysis and series inversions, Springer-Verlag, 1985.

\bibitem{Chen93}
W. Y. C. Chen, Context-free grammars, differential operators
and formal power series, Theoret. Comput. Sci. {\bf 117} (1993), 113-129.

\bibitem{Dumont}
D. Dumont and A. Ramamonjisoa, Grammaire de
Ramanujan et Arbres de Cayley, Electron. J. Combin. {\bf 3} No. 2 (1996),
R17.

\bibitem{Shor95}
P. Shor,  A new proof of Cayley's formula for counting labeled trees,  
J. Combin. Theory. Ser. A {\bf 71} (1995), 154-158.

\bibitem{Zeng}
J. Zeng, A Ramanujan sequence that refines the Cayley formula for trees, 
Ramanujan J. {\bf 3} (1999), 45--54.
\end{thebibliography}
\end{document}